%Paper: alg-geom/9405012
%From: laszlo@ceremab.u-bordeaux.fr (Yves Laszlo)
%Date: Tue, 24 May 94 11:45:00 +0200

\magnification=1250
\hsize=13.5cm
\vsize=20cm

\parindent=12pt   \parskip=0pt
\pageno=1

\hoffset=0mm
\voffset=0cm

\ifnum\mag=\magstep1
\hoffset=0cm   % offset horizontal en \magnification=1200
\voffset=-0,5cm   % offset horizontal en \magnification=1200
\fi

\pretolerance=500 \tolerance=1000  \brokenpenalty=5000

\catcode`\@=11

\font\eightrm=cmr8         \font\eighti=cmmi8
\font\eightsy=cmsy8        \font\eightbf=cmbx8
\font\eighttt=cmtt8        \font\eightit=cmti8
\font\eightsl=cmsl8        \font\sixrm=cmr6
\font\sixi=cmmi6           \font\sixsy=cmsy6
\font\sixbf=cmbx6

\font\tengoth=eufm10       \font\tenbboard=msbm10
\font\eightgoth=eufm8      \font\eightbboard=msbm8
\font\sevengoth=eufm7      \font\sevenbboard=msbm7
\font\sixgoth=eufm6        \font\fivegoth=eufm5

\font\tencyr=wncyr10       %\font\tenbcyr=wncyb10
\font\eightcyr=wncyr8      %\font\eightbcyr=wncyb8
\font\sevencyr=wncyr7      %\font\sevenbcyr=wncyb7
\font\sixcyr=wncyr6

\font\tengo=eufm10 \font\sevengo=eufm7 \font\fivego=eufm5
\newfam\gofam \textfont\gofam=\tengo \scriptfont\gofam=\sevengo
\scriptscriptfont\gofam=\fivego 
% Pour que les accents se placent correctement en mode math en corps 8 et 6

\skewchar\eighti='177 \skewchar\sixi='177
\skewchar\eightsy='60 \skewchar\sixsy='60

\newfam\gothfam           \newfam\bboardfam
\newfam\cyrfam

\def\tenpoint{%
  \textfont0=\tenrm \scriptfont0=\sevenrm \scriptscriptfont0=\fiverm
  \def\rm{\fam\z@\tenrm}%
  \textfont1=\teni  \scriptfont1=\seveni  \scriptscriptfont1=\fivei
  \def\oldstyle{\fam\@ne\teni}\let\old=\oldstyle
  \textfont2=\tensy \scriptfont2=\sevensy \scriptscriptfont2=\fivesy
  \textfont\gothfam=\tengoth \scriptfont\gothfam=\sevengoth
  \scriptscriptfont\gothfam=\fivegoth
  \def\goth{\fam\gothfam\tengoth}%
  \textfont\bboardfam=\tenbboard \scriptfont\bboardfam=\sevenbboard
  \scriptscriptfont\bboardfam=\sevenbboard
  \def\bb{\fam\bboardfam\tenbboard}%
 \textfont\cyrfam=\tencyr \scriptfont\cyrfam=\sevencyr
  \scriptscriptfont\cyrfam=\sixcyr
  \def\cyr{\fam\cyrfam\tencyr}%
  \textfont\itfam=\tenit
  \def\it{\fam\itfam\tenit}%
  \textfont\slfam=\tensl
  \def\sl{\fam\slfam\tensl}%
  \textfont\bffam=\tenbf \scriptfont\bffam=\sevenbf
  \scriptscriptfont\bffam=\fivebf
  \def\bf{\fam\bffam\tenbf}%
  \textfont\ttfam=\tentt
  \def\tt{\fam\ttfam\tentt}%
  \abovedisplayskip=12pt plus 3pt minus 9pt
  \belowdisplayskip=\abovedisplayskip
  \abovedisplayshortskip=0pt plus 3pt
  \belowdisplayshortskip=4pt plus 3pt
  \smallskipamount=3pt plus 1pt minus 1pt
  \medskipamount=6pt plus 2pt minus 2pt
  \bigskipamount=12pt plus 4pt minus 4pt
  \normalbaselineskip=12pt
  \setbox\strutbox=\hbox{\vrule height8.5pt depth3.5pt width0pt}%
  \let\bigf@nt=\tenrm       \let\smallf@nt=\sevenrm
  \normalbaselines\rm}

\def\eightpoint{%
  \textfont0=\eightrm \scriptfont0=\sixrm \scriptscriptfont0=\fiverm
  \def\rm{\fam\z@\eightrm}%
  \textfont1=\eighti  \scriptfont1=\sixi  \scriptscriptfont1=\fivei
  \def\oldstyle{\fam\@ne\eighti}\let\old=\oldstyle
  \textfont2=\eightsy \scriptfont2=\sixsy \scriptscriptfont2=\fivesy
  \textfont\gothfam=\eightgoth \scriptfont\gothfam=\sixgoth
  \scriptscriptfont\gothfam=\fivegoth
  \def\goth{\fam\gothfam\eightgoth}%
  \textfont\cyrfam=\eightcyr \scriptfont\cyrfam=\sixcyr
  \scriptscriptfont\cyrfam=\sixcyr
  \def\cyr{\fam\cyrfam\eightcyr}%
  \textfont\bboardfam=\eightbboard \scriptfont\bboardfam=\sevenbboard
  \scriptscriptfont\bboardfam=\sevenbboard
  \def\bb{\fam\bboardfam}%
  \textfont\itfam=\eightit
  \def\it{\fam\itfam\eightit}%
  \textfont\slfam=\eightsl
  \def\sl{\fam\slfam\eightsl}%
  \textfont\bffam=\eightbf \scriptfont\bffam=\sixbf
  \scriptscriptfont\bffam=\fivebf
  \def\bf{\fam\bffam\eightbf}%
  \textfont\ttfam=\eighttt
  \def\tt{\fam\ttfam\eighttt}%
  \abovedisplayskip=9pt plus 3pt minus 9pt
  \belowdisplayskip=\abovedisplayskip
  \abovedisplayshortskip=0pt plus 3pt
  \belowdisplayshortskip=3pt plus 3pt
  \smallskipamount=2pt plus 1pt minus 1pt
  \medskipamount=4pt plus 2pt minus 1pt
  \bigskipamount=9pt plus 3pt minus 3pt
  \normalbaselineskip=9pt
  \setbox\strutbox=\hbox{\vrule height7pt depth2pt width0pt}%
  \let\bigf@nt=\eightrm     \let\smallf@nt=\sixrm
  \normalbaselines\rm}

\tenpoint

\def\pc#1{\bigf@nt#1\smallf@nt}         \def\pd#1 {{\pc#1} }

\catcode`\;=\active
\def;{\relax\ifhmode\ifdim\lastskip>\z@\unskip\fi
\kern\fontdimen2  -1.2 \fontdimen3 \string;}

\catcode`\:=\active
\def:{\relax\ifhmode\ifdim\lastskip>\z@\unskip\fi\penalty\@M\ \fi\string:}

\catcode`\!=\active
\def!{\relax\ifhmode\ifdim\lastskip>\z@
\unskip\fi\kern\fontdimen2  -1.1 \fontdimen3 \string!}

\catcode`\?=\active
\def?{\relax\ifhmode\ifdim\lastskip>\z@
\unskip\fi\kern\fontdimen2  -1.1 \fontdimen3 \string?}

\def\^#1{\if#1i{\accent"5E\i}\else{\accent"5E #1}\fi}
\def\"#1{\if#1i{\accent"7F\i}\else{\accent"7F #1}\fi}

\frenchspacing

\newtoks\auteurcourant      \auteurcourant={\hfil}
\newtoks\titrecourant       \titrecourant={\hfil}

\newtoks\hautpagetitre      \hautpagetitre={\hfil}
\newtoks\baspagetitre       \baspagetitre={\hfil}

\newtoks\hautpagegauche
\hautpagegauche={\eightpoint\rlap{\folio}\hfil\the\auteurcourant\hfil}
\newtoks\hautpagedroite
\hautpagedroite={\eightpoint\hfil\the\titrecourant\hfil\llap{\folio}}

\newtoks\baspagegauche      \baspagegauche={\hfil}
\newtoks\baspagedroite      \baspagedroite={\hfil}

\newif\ifpagetitre          \pagetitretrue

% \nopagenumbers : c'est un peu violent, mais ça marche. Alors ...

\headline={\ifpagetitre\the\hautpagetitre
\else\ifodd\pageno\the\hautpagedroite\else\the\hautpagegauche\fi\fi}

\footline={\ifpagetitre\the\baspagetitre\else
\ifodd\pageno\the\baspagedroite\else\the\baspagegauche\fi\fi
\global\pagetitrefalse}

% Redefinition de \raggedbottom pour avoir plus de mou en bas de page
% (necesssaire quand il y a beaucoup de grumeaux, des grosses
% formules centrees et pas beaucoup de texte entre)

\def\raggedbottom{\topskip 10pt plus 36pt\r@ggedbottomtrue}

\def\pointir{\unskip . --- \ignorespaces}

\def\Bigbreak{\vskip-\lastskip\bigbreak}
\def\Medbreak{\vskip-\lastskip\medbreak}

\def\ctexte#1\endctexte{%
  \hbox{$\vcenter{\halign{\hfill##\hfill\crcr#1\crcr}}$}}

\long\def\ctitre#1\endctitre{%
    \ifdim\lastskip<24pt\vskip-\lastskip\bigbreak\bigbreak\fi
  		\vbox{\parindent=0pt\leftskip=0pt plus 1fill
          \rightskip=\leftskip
          \parfillskip=0pt\bf#1\par}
    \bigskip\nobreak}

\long\def\section#1\endsection{%
\vskip 0pt plus 3\normalbaselineskip
\penalty-250
\vskip 0pt plus -3\normalbaselineskip
\Bigbreak
\message{[section \string: #1]}{\bf#1\unskip}\pointir}

\long\def\sectiona#1\endsection{%
\vskip 0pt plus 3\normalbaselineskip
\penalty-250
\vskip 0pt plus -3\normalbaselineskip
\Bigbreak
\message{[sectiona \string: #1]}%
{\bf#1}\medskip\nobreak}

\long\def\subsection#1\endsubsection{%
\Medbreak
{\it#1\unskip}\pointir}

\long\def\subsectiona#1\endsubsection{%
\Medbreak
{\it#1}\par\nobreak}

\def\rem#1\endrem{%
\Medbreak
{\it#1\unskip} : }

\def\remp#1\endrem{%
\Medbreak
{\pc #1\unskip}\pointir}

\def\rema#1\endrem{%
\Medbreak
{\it #1}\par\nobreak}

\def\newparwithcolon#1\endnewparwithcolon{
\Medbreak
{#1\unskip} : }

\def\newparwithpointir#1\endnewparwithpointir{
\Medbreak
{#1\unskip}\pointir}

\def\newpara#1\endnewpar{
\Medbreak
{#1\unskip}\smallskip\nobreak}

\long\def\th#1 #2\enonce#3\endth{%
   \Medbreak
   {\pc#1} {#2\unskip}\pointir{\it #3}\medskip}

\long\def\tha#1 #2\enonce#3\endth{%
   \Medbreak
   {\pc#1} {#2\unskip}\par\nobreak{\it #3}\medskip}

\long\def\Th#1 #2 #3\enonce#4\endth{%
   \Medbreak
   #1 {\pc#2} {#3\unskip}\pointir{\it #4}\medskip}

\long\def\Tha#1 #2 #3\enonce#4\endth{%
   \Medbreak
   #1 {\pc#2} #3\par\nobreak{\it #4}\medskip}

\def\decale#1{\smallbreak\hskip 28pt\llap{#1}\kern 5pt}
\def\decaledecale#1{\smallbreak\hskip 34pt\llap{#1}\kern 5pt}
\def\puce{\smallbreak\hskip 6pt{$\scriptstyle\bullet$}\kern 5pt}

\def\displaylinesno#1{\displ@y\halign{
\hbox to\displaywidth{$\@lign\hfil\displaystyle##\hfil$}&
\llap{$##$}\crcr#1\crcr}}

\def\ldisplaylinesno#1{\displ@y\halign{
\hbox to\displaywidth{$\@lign\hfil\displaystyle##\hfil$}&
\kern-\displaywidth\rlap{$##$}\tabskip\displaywidth\crcr#1\crcr}}

\def\eqalign#1{\null\,\vcenter{\openup\jot\m@th\ialign{
\strut\hfil$\displaystyle{##}$&$\displaystyle{{}##}$\hfil
&&\quad\strut\hfil$\displaystyle{##}$&$\displaystyle{{}##}$\hfil
\crcr#1\crcr}}\,}

\def\system#1{\left\{\null\,\vcenter{\openup1\jot\m@th
\ialign{\strut$##$&\hfil$##$&$##$\hfil&&
        \enskip$##$\enskip&\hfil$##$&$##$\hfil\crcr#1\crcr}}\right.}

\let\@ldmessage=\message

\def\message#1{{\def\pc{\string\pc\space}%
                \def\'{\string'}\def\`{\string`}%
                \def\^{\string^}\def\"{\string"}%
                \@ldmessage{#1}}}

\def\up#1{\raise 1ex\hbox{\smallf@nt#1}}

\def\qed{\raise -2pt\hbox{\vrule\vbox to 10pt{\hrule width 4pt
                 \vfill\hrule}\vrule}}

\def\cqfd{\unskip\penalty 500\quad\qed\medbreak}

\def\virg{\raise .4ex\hbox{,}}   % virgule après une fraction

 % point-virgule de ponctuation en maths

\def\build#1_#2^#3{\mathrel{
\mathop{\kern 0pt#1}\limits_{#2}^{#3}}}

\def\boxit#1#2{%
\setbox1=\hbox{\kern#1{#2}\kern#1}%
\dimen1=\ht1 \advance\dimen1 by #1 \dimen2=\dp1 \advance\dimen2 by #1
\setbox1=\hbox{\vrule height\dimen1 depth\dimen2\box1\vrule}%
\setbox1=\vbox{\hrule\box1\hrule}%
\advance\dimen1 by .4pt \ht1=\dimen1
\advance\dimen2 by .4pt \dp1=\dimen2  \box1\relax}

\def\fhd#1#2{\smash{\mathop{\hbox to 12mm{\rightarrowfill}}
\limits^{\scriptstyle#1}_{\scriptstyle#2}}}
\def\fhg#1#2{\smash{\mathop{\hbox to 12mm{\leftarrowfill}}
\limits^{\scriptstyle#1}_{\scriptstyle#2}}}

\frenchspacing
\parindent=0mm
\baspagegauche={\centerline{\tenbf\folio}}
\baspagedroite={\centerline{\tenbf\folio}}
\hautpagegauche={\hfil}
\hautpagedroite={\hfil}
\def\ind{\hskip 1cm\relax}
\def\Z{\hbox{\bb Z}}
\def\ra{\rightarrow}

\catcode`\@=12

\showboxbreadth=-1  \showboxdepth=-1

\def\Grille{\setbox13=\vbox to 5\unitlength{\hrule width 109mm\vfill}
\setbox13=\vbox to 65mm{\offinterlineskip\leaders\copy13\vfill\kern-1pt\hrule}
\setbox14=\hbox to 5\unitlength{\vrule height 65mm\hfill}
\setbox14=\hbox to 109mm{\leaders\copy14\hfill\kern-2mm\vrule height 65mm}
\ht14=0pt\dp14=0pt\wd14=0pt \setbox13=\vbox to
0pt{\vss\box13\offinterlineskip\box14} \wd13=0pt\box13}

\def\fleche(#1,#2)\dir(#3,#4)\long#5{%
\noalign{\leftput(#1,#2){\vector(#3,#4){#5}}}}

\def\ligne(#1,#2)\dir(#3,#4)\long#5{%
\noalign{\leftput(#1,#2){\lline(#3,#4){#5}}}}

\def\put(#1,#2)#3{\noalign{\setbox1=\hbox{%
    \kern #1\unitlength
    \raise #2\unitlength\hbox{$#3$}}%
    \ht1=0pt \wd1=0pt \dp1=0pt\box1}}

\def\hfl#1#2#3{\smash{\mathop{\hbox to#3{\rightarrowfill}}\limits
^{\scriptstyle#1}_{\scriptstyle#2}}}

\def\gfl#1#2#3{\smash{\mathop{\hbox to#3{\leftarrowfill}}\limits
^{\scriptstyle#1}_{\scriptstyle#2}}}

 \message{`lline' & `vector' macros from LaTeX}
 \catcode`@=11
\def\{{\relax\ifmmode\lbrace\else$\lbrace$\fi}
\def\}{\relax\ifmmode\rbrace\else$\rbrace$\fi}
\def\newcount{\alloc@0\count\countdef\insc@unt}
\def\newdimen{\alloc@1\dimen\dimendef\insc@unt}
\def\newwrite{\alloc@7\write\chardef\sixt@@n}

\newwrite\@unused
\newcount\@tempcnta
\newcount\@tempcntb
\newdimen\@tempdima
\newdimen\@tempdimb
\newbox\@tempboxa

\def\@spaces{\space\space\space\space}
\def\@whilenoop#1{}
\def\@whiledim#1\do #2{\ifdim #1\relax#2\@iwhiledim{#1\relax#2}\fi}
\def\@iwhiledim#1{\ifdim #1\let\@nextwhile=\@iwhiledim
        \else\let\@nextwhile=\@whilenoop\fi\@nextwhile{#1}}
\def\@badlinearg{\@latexerr{Bad \string\line\space or \string\vector
   \space argument}}
\def\@latexerr#1#2{\begingroup
\edef\@tempc{#2}\expandafter\errhelp\expandafter{\@tempc}%
%% error help message pieces.
\def\@eha{Your command was ignored.
^^JType \space I <command> <return> \space to replace it
  with another command,^^Jor \space <return> \space to continue without
it.}
\def\@ehb{You've lost some text. \space \@ehc}
\def\@ehc{Try typing \space <return>
  \space to proceed.^^JIf that doesn't work, type \space X <return> \space to
  quit.}
\def\@ehd{You're in trouble here.  \space\@ehc}

\typeout{LaTeX error. \space See LaTeX manual for explanation.^^J
 \space\@spaces\@spaces\@spaces Type \space H <return> \space for
 immediate help.}\errmessage{#1}\endgroup}
\def\typeout#1{{\let\protect\string\immediate\write\@unused{#1}}}

% line & circle fonts
\font\tenln    = line10
\font\tenlnw   = linew10
%\font\tencirc  = circle10
%\font\tencircw = circlew10

\newdimen\@wholewidth
\newdimen\@halfwidth
\newdimen\unitlength

\unitlength =1pt

%\newbox\@picbox
%\newdimen\@picht

\def\thinlines{\let\@linefnt\tenln \let\@circlefnt\tencirc
  \@wholewidth\fontdimen8\tenln \@halfwidth .5\@wholewidth}
\def\thicklines{\let\@linefnt\tenlnw \let\@circlefnt\tencircw
  \@wholewidth\fontdimen8\tenlnw \@halfwidth .5\@wholewidth}

\def\linethickness#1{\@wholewidth #1\relax \@halfwidth .5\@wholewidth}

\newif\if@negarg

\def\lline(#1,#2)#3{\@xarg #1\relax \@yarg #2\relax
\@linelen=#3\unitlength
\ifnum\@xarg =0 \@vline
  \else \ifnum\@yarg =0 \@hline \else \@sline\fi
\fi}

\def\@sline{\ifnum\@xarg< 0 \@negargtrue \@xarg -\@xarg \@yyarg -\@yarg
  \else \@negargfalse \@yyarg \@yarg \fi
\ifnum \@yyarg >0 \@tempcnta\@yyarg \else \@tempcnta -\@yyarg \fi
\ifnum\@tempcnta>6 \@badlinearg\@tempcnta0 \fi
\setbox\@linechar\hbox{\@linefnt\@getlinechar(\@xarg,\@yyarg)}%
\ifnum \@yarg >0 \let\@upordown\raise \@clnht\z@
   \else\let\@upordown\lower \@clnht \ht\@linechar\fi
\@clnwd=\wd\@linechar
\if@negarg \hskip -\wd\@linechar \def\@tempa{\hskip -2\wd\@linechar}\else
     \let\@tempa\relax \fi
\@whiledim \@clnwd <\@linelen \do
  {\@upordown\@clnht\copy\@linechar
   \@tempa
   \advance\@clnht \ht\@linechar
   \advance\@clnwd \wd\@linechar}%
\advance\@clnht -\ht\@linechar
\advance\@clnwd -\wd\@linechar
\@tempdima\@linelen\advance\@tempdima -\@clnwd
\@tempdimb\@tempdima\advance\@tempdimb -\wd\@linechar
\if@negarg \hskip -\@tempdimb \else \hskip \@tempdimb \fi
\multiply\@tempdima \@m
\@tempcnta \@tempdima \@tempdima \wd\@linechar \divide\@tempcnta \@tempdima
\@tempdima \ht\@linechar \multiply\@tempdima \@tempcnta
\divide\@tempdima \@m
\advance\@clnht \@tempdima
\ifdim \@linelen <\wd\@linechar
   \hskip \wd\@linechar
  \else\@upordown\@clnht\copy\@linechar\fi}

\def\@hline{\ifnum \@xarg <0 \hskip -\@linelen \fi
\vrule height \@halfwidth depth \@halfwidth width \@linelen
\ifnum \@xarg <0 \hskip -\@linelen \fi}

\def\@getlinechar(#1,#2){\@tempcnta#1\relax\multiply\@tempcnta 8
\advance\@tempcnta -9 \ifnum #2>0 \advance\@tempcnta #2\relax\else
\advance\@tempcnta -#2\relax\advance\@tempcnta 64 \fi
\char\@tempcnta}

\def\vector(#1,#2)#3{\@xarg #1\relax \@yarg #2\relax
\@linelen=#3\unitlength
\ifnum\@xarg =0 \@vvector
  \else \ifnum\@yarg =0 \@hvector \else \@svector\fi
\fi}

\def\@hvector{\@hline\hbox to 0pt{\@linefnt
\ifnum \@xarg <0 \@getlarrow(1,0)\hss\else
    \hss\@getrarrow(1,0)\fi}}

\def\@vvector{\ifnum \@yarg <0 \@downvector \else \@upvector \fi}

\def\@svector{\@sline
\@tempcnta\@yarg \ifnum\@tempcnta <0 \@tempcnta=-\@tempcnta\fi
\ifnum\@tempcnta <5
  \hskip -\wd\@linechar
  \@upordown\@clnht \hbox{\@linefnt  \if@negarg
  \@getlarrow(\@xarg,\@yyarg) \else \@getrarrow(\@xarg,\@yyarg) \fi}%
\else\@badlinearg\fi}

\def\@getlarrow(#1,#2){\ifnum #2 =\z@ \@tempcnta='33\else
\@tempcnta=#1\relax\multiply\@tempcnta \sixt@@n \advance\@tempcnta
-9 \@tempcntb=#2\relax\multiply\@tempcntb \tw@
\ifnum \@tempcntb >0 \advance\@tempcnta \@tempcntb\relax
\else\advance\@tempcnta -\@tempcntb\advance\@tempcnta 64
\fi\fi\char\@tempcnta}

\def\@getrarrow(#1,#2){\@tempcntb=#2\relax
\ifnum\@tempcntb < 0 \@tempcntb=-\@tempcntb\relax\fi
\ifcase \@tempcntb\relax \@tempcnta='55 \or
\ifnum #1<3 \@tempcnta=#1\relax\multiply\@tempcnta
24 \advance\@tempcnta -6 \else \ifnum #1=3 \@tempcnta=49
\else\@tempcnta=58 \fi\fi\or
\ifnum #1<3 \@tempcnta=#1\relax\multiply\@tempcnta
24 \advance\@tempcnta -3 \else \@tempcnta=51\fi\or
\@tempcnta=#1\relax\multiply\@tempcnta
\sixt@@n \advance\@tempcnta -\tw@ \else
\@tempcnta=#1\relax\multiply\@tempcnta
\sixt@@n \advance\@tempcnta 7 \fi\ifnum #2<0 \advance\@tempcnta 64 \fi
\char\@tempcnta}

\def\@vline{\ifnum \@yarg <0 \@downline \else \@upline\fi}

\def\@upline{\hbox to \z@{\hskip -\@halfwidth \vrule
  width \@wholewidth height \@linelen depth \z@\hss}}

\def\@downline{\hbox to \z@{\hskip -\@halfwidth \vrule
  width \@wholewidth height \z@ depth \@linelen \hss}}

\def\@upvector{\@upline\setbox\@tempboxa\hbox{\@linefnt\char'66}\raise
     \@linelen \hbox to\z@{\lower \ht\@tempboxa\box\@tempboxa\hss}}

\def\@downvector{\@downline\lower \@linelen
      \hbox to \z@{\@linefnt\char'77\hss}}

\thinlines

\newcount\@xarg
\newcount\@yarg
\newcount\@yyarg
\newcount\@multicnt
\newdimen\@xdim
\newdimen\@ydim
\newbox\@linechar
\newdimen\@linelen
\newdimen\@clnwd
\newdimen\@clnht
\newdimen\@dashdim
\newbox\@dashbox
\newcount\@dashcnt
 \catcode`@=12

\newbox\tbox
\newbox\tboxa

\def\leftzer#1{\setbox\tbox=\hbox to 0pt{#1\hss}%
     \ht\tbox=0pt \dp\tbox=0pt \box\tbox}

\def\rightzer#1{\setbox\tbox=\hbox to 0pt{\hss #1}%
     \ht\tbox=0pt \dp\tbox=0pt \box\tbox}

\def\centerzer#1{\setbox\tbox=\hbox to 0pt{\hss #1\hss}%
     \ht\tbox=0pt \dp\tbox=0pt \box\tbox}

\def\image(#1,#2)#3{\vbox to #1{\offinterlineskip
    \vss #3 \vskip #2}}

\def\leftput(#1,#2)#3{\setbox\tboxa=\hbox{%
    \kern #1\unitlength
    \raise #2\unitlength\hbox{\leftzer{#3}}}%
    \ht\tboxa=0pt \wd\tboxa=0pt \dp\tboxa=0pt\box\tboxa}

\def\rightput(#1,#2)#3{\setbox\tboxa=\hbox{%
    \kern #1\unitlength
    \raise #2\unitlength\hbox{\rightzer{#3}}}%
    \ht\tboxa=0pt \wd\tboxa=0pt \dp\tboxa=0pt\box\tboxa}

\def\centerput(#1,#2)#3{\setbox\tboxa=\hbox{%
    \kern #1\unitlength
    \raise #2\unitlength\hbox{\centerzer{#3}}}%
    \ht\tboxa=0pt \wd\tboxa=0pt \dp\tboxa=0pt\box\tboxa}

\unitlength=1mm

% MAJMATH.TEX

\mathcode`A="7041 \mathcode`B="7042 \mathcode`C="7043 \mathcode`D="7044
\mathcode`E="7045 \mathcode`F="7046 \mathcode`G="7047 \mathcode`H="7048
\mathcode`I="7049 \mathcode`J="704A \mathcode`K="704B \mathcode`L="704C
\mathcode`M="704D \mathcode`N="704E \mathcode`O="704F \mathcode`P="7050
\mathcode`Q="7051 \mathcode`R="7052 \mathcode`S="7053 \mathcode`T="7054
\mathcode`U="7055 \mathcode`V="7056 \mathcode`W="7057 \mathcode`X="7058
\mathcode`Y="7059 \mathcode`Z="705A

\def\spacedmath#1{\def\packedmath##1${\bgroup\mathsurround=0pt ##1\egroup$}%
\mathsurround#1 \everymath={\packedmath}\everydisplay={\mathsurround=0pt }}

\def\nospacedmath{\mathsurround=0pt \everymath={}\everydisplay={} }

\def\ind{\hskip 1cm\relax}

\def\dra{\rightarrow\kern -3mm\rightarrow}

\def\limproj{\mathop{\oalign{lim\cr\hidewidth$\longleftarrow$\hidewidth\cr}}}

\def\Z{\hbox{\bb Z}}

\def\fhd#1#2{\nospacedmath\smash{\mathop{\hbox to 12mm{\rightarrowfill}}
\limits^{\scriptstyle#1}_{\scriptstyle#2}}}
\def\fhg#1#2{\nospacedmath\smash{\mathop{\hbox to 12mm{\leftarrowfill}}
\limits^{\scriptstyle#1}_{\scriptstyle#2}}}

\def\mv=msam10\def\dra{{\mv\char 16}}

\spacedmath{2pt}

\parindent=0mm
\def\rond{\kern 1pt{\scriptstyle\circ}\kern 1pt}
\def\epi{\rightarrow \kern -3mm\rightarrow }
\def\longepi{\longrightarrow \kern -5mm\longrightarrow }
\def\mono{\lhook\joinrel\mathrel{\longrightarrow}}

\def\Ext{\mathop{\rm Ext}\nolimits}

\def\Aut{\mathop{\rm Aut}\nolimits}

\def\Ker{\mathop{\rm Ker}\nolimits}

\def\det{\mathop{\rm det}\nolimits}
\def\Gl{\mathop{\rm Gl}\nolimits}
\def\O{{\cal O}}
\def\mult{\mathop{\rm mult}\nolimits}
\def\Sl{\mathop{\rm Sl}\nolimits}
\def\Pic{\mathop{\rm Pic}\nolimits}

\def\dim{\mathop{\rm dim}\nolimits}

\def\Sp{\mathop{\rm Spec\,}\nolimits}
\def\iso{\mathrel{\mathop{\kern 0pt\longrightarrow }\limits^{\sim}}}
\def\div{\mathop{\rm div\,}\nolimits}

\baselineskip 15pt

\overfullrule=0pt

\catcode`\@=12

\showboxbreadth=-1  \showboxdepth=-1

\def\bl{\bigl}\def\Bl{\Bigl}\def\br{\bigr}
\def\Br{\Bigr}

\def\P{{\bf P}}
\def\S{{\cal S}}
\def\int{{\kern 1truemm \hbox{\vrule height .4pt width 7pt}\hbox{\vrule height
7pt width .4pt}\kern
1truemm}}
\def\ind{\par\hskip 1cm\relax}

\unitlength=1mm
\vsize = 25truecm
\hsize = 16truecm
\hoffset = -.15truecm
\voffset = -.5truecm

\parskip=0pt
\pageno=1

\pretolerance=500 \tolerance=1000  \brokenpenalty=5000
\def\cqfd{\kern 2truemm\unskip\penalty 500\vrule
height 4pt depth 0pt width 4pt\medbreak}

\def\D{{{\cal O}(dx)}}

\def\m{{SU(r,d)}}

\parindent=0cm
\def\K{{\tilde K}}

\centerline{{\bf Local structure of the moduli space of vector bundles over
curves}}
\medskip

\centerline{Yves {\pc LASZLO}
\footnote{(*)}{\sevenrm The author was partially
supported by the European Science Project ``Geometry of Algebraic Varieties",
Contract no. SCI-0398-C(A).}} \vskip1.5cm
\parindent=0cm

\section 0. Introduction
\endsection
\ind Let $X$ be a smooth, projective and connected curve (over an algebraically
closed field of
characteristic zero) of genus $g(X)\geq 2$. Let $x$ be a (closed) point
of $X$ and $\m$ the moduli space of semi-stable vector bundles on $X$ of rank
$r\geq 2$ and
determinant $\D$. As usual, the geometric points of $\m$ correspond to
$S$-equivalence classes
$[E]$ where $E$ is a semi-stable rank $r$ bundle of determinant $\D$ (another
semi-stable bundle
$F$ is said to be $S$-equivalent to $F$ if the graded objects $gr(E)$ and
$gr(F)$ are
isomorphic).

\ind The singular locus of $\m$ consists exactly of the non stable points
(except if $r=g(X)=2$
and $d=0$. In this case, $\m=\P^3$ [N-R1]). In particular,
except in the exceptional case above, $\m$ is smooth if and only if $r$ and $d$
are
relatively prime. General facts about the action of reductive groups ensure
that $\m$ is
Cohen-Macaulay [E-H], normal and that the singularities are rational [B]. The
principal aim of
this paper is to give additional information about the singularities,
essentially the
description of the completion  of the
local ring at a non smooth point of $\m$ and to compute the multiplicity and
the tangent cones at
those singular points $[E]$ which are not too bad, {\it i.e.} the corresponding
graded object
$gr(E)$ of $[E]$ has only two non isomorphic stable summands  (or equivalently
$\Aut\bl(gr(E)\br)=G_m\times G_m$). Further, we give a complete description in
the rank 2 case
(proposition II.2, corollary II.3 and theorem III.4).

\ind As a corollary, we get the local form of the so called
Coble quartic and prove that the Kummer variety of the  Jacobian of a genus 3
non hyperelliptic
curve is {\it schematically} defined by 8 cubics, the partials derivatives of
the Coble
quartic (theorem III.6).

\ind One could also give partial information at least if $\Aut\bl(gr(E)\br)$ is
a torus, or
by using results of [P], if $\Aut\bl(gr(E)\br)=\Gl_r(k)$ (the latter case
essentialy means that
$gr(E)$ is  the trivial bundle). But it seems to be difficult and somewhat
messy to calculate
for instance the multiplicity. In the remaining part of the paper, we compute
the multiplicity of
a generalized theta divisor of $SU_X(2,\O)$ at a point $[L\oplus L^{\vee}]$,
where $L^{2}\not=\O$.
In fact, this computation could be done with only minor changes for a point
$[E]$ of any rank
whith $\det(gr(E))=\O$ and $\Aut\bl(gr(E)\br)=G_m\times G_m$.

\ind Let us also mention that similar results could be obtained exactly in the
same way for
certain surfaces. But, all the future applications that we have in mind as well
as the
applications that we have in our hand are for curves. Therefore, we have
restricted ourselves
to the case of curves.
\medskip

{\sevenrm I would like to
thank {\pc A.} {\pc BEAUVILLE} for usefull discussions and {\pc M.} {\pc VAN}
{\pc DER} {\pc
PUT} for helping me to correct the english of a preliminary version.}
\medskip
{\bf Notations and conventions.} All the schemes are
of finite type over $k$, by point we mean closed point. If $(X_i)_{1\leq N}$
(resp. $(n_i)_{1\leq
N}$) are  indeterminates (resp. non negative integers), let me denote by
$\underline X$,
$\underline n$ and $X^{\underline n}$ the $N-$tuple $\underline X=(X_i)$, the
multi-index
$\underline n=(n_i)$  and the product $X^{\underline n}=\prod_{i=1}^Nx_i^{n_i}$
repectively. For
$V$  a finite dimensional vector space with dual $V^\vee$, the ring $k[V]$
(resp. $k[[V]]$) is the
polynomial ring ${\rm Sym}V^\vee$ (resp. its completion at the origin). The
scheme $\P(V)$ is
the projective space ${\rm Proj}({\rm Sym}V^\vee)$ of lines of $V$ and $V$ will
also denote the
pointed affine space $\Sp({\rm Sym}V^\vee)=\Sp k[V]$ (notice that $k[V]$ is the
coordinate ring of
$V$). Let $(E_i)$ be a set of vector bundles over $X$, the kernel of the trace
map
$$\Ker(\oplus\Ext^1(E_i,E_i)\fhd{\oplus Tr_i}{}H^1(X,\O))$$ will be denoted by
$$\bl(\oplus\Ext^1(gr_i,gr_i)\br)_0.$$
Finally, $E$ will always denote a rank $r$ semi-stable bundle on $X$ of degree
$d$. Let $Fil$ be a
strictly increasing Jordan-H\"older filtration of $E$ by stable bundles
$Fil^i(E)$ (with  slope
${d\over r}$.) $$0=Fil^0(E)\subset Fil^1(E)\ldots\subset Fil^N(E)=E.$$ Then,
the graded object
$$\oplus gr_i$$ with
$$gr_i=Fil^i(E)/Fil^{i-1}(E)$$
is well defined (up to isomorphism!) and is denoted by $gr$.

\section I. Local structure of $\m$ and classical invariant theory
\endsection

\ind It is well known (see [S] for instance) that the key ingredient to analyse
the local
structure of $\m$ is the \'etale slice theorem of Luna. Let us recall this
analysis.

Take $n$ big enough such that $E(nx)$ is globally generated and has no $H^1$
for every semi-stable
rank $r$ vector bundle $E$ of degree $d$ (every $n$ such that $rn+d>r(2g-1)$
has this property). Let
$\chi=\chi(E(nx))$ be the correponding Euler-characteristic. In Grothendieck's
scheme ${\cal
Q}uot$ which parametrizes quotients
$$\O(-nx)^{\oplus \chi}\epi E,$$
let $\cal Q$ be the open set whose closed
points correspond to such quotients with the following properties:

(i) $E$ is semi-stable of rank $r$ and degree $d$.

(ii) $H^1(X,E(nx))=0$ and the natural map $H^0(X,\O^{\oplus \chi})\ra
H^0(X,E(nx))$ is onto.
\medskip
Let $\cal E$ be the universal quotient bundle on $\cal Q$.The scheme $\cal Q$
is smooth and
the semi-simple group $G=PGL_{\chi}$ acts
on it.The moduli space $\m$ is the GIT quotient of ${\cal Q}//G$.
\medskip
\ind Let $q=[\O(-nx)^{\oplus \chi}\epi E]$ a point of $\cal Q$ and $gr$ the
corresponding graded
object. By the very definition of semi-stability, there exists a $G$-stable
open affine
neighborhood $\Omega$ and the fibre of $\Omega\ra\Omega/G$ at $[E]$ contains a
unique closed
orbit $G(q)$. This orbit is either characterized as being of minimal dimension,
or as having an isotropy group $G_q=\Aut(E)$ of maximal dimension. Note
that the scalar matrices act trivially.  Inside the isotropy groups of
the elements of the S-equivalence class of $E$, the isotropy group $\Aut(gr)$
is of maximal
dimension, and therefore the corresponding orbit is closed. The closedness of
the orbit
corresponding to $gr$ allows us to use the Luna \'etale slice theorem
[Lu] which gives precisely the:

\th THEOREM 1 (Luna)
\enonce There exists a
closed subscheme $V$ of $\Omega$ such that:

(i) $V$ is stable under the isotropy group $G_q$ and the tangent space $T_qV$
is
$G_q-$\nobreak isomorphic to $T_q\Omega/T_qG(q)=\Ext_0^1(gr,gr)$.

(ii) The $G$-morphism
$$G\times_{G_x}V\ra G.V$$
is \'etale and onto the $G$-saturated {\it open} subset $G(V)$ of $\m$.

(iii) The induced morphism $V/G_q\ra \Omega/G$ is \'etale at the image $\bar q$
of $q$ in
$V/G_q$.
\endth
\ind Let $G$ be a reductive group acting on a $k$-algebra of finite type $R$
and
$\goth m$ an invariant ideal. Let $R_0$ be the ring of invariant and
${\goth m}_0={\goth m}\cap R_0$ the corresponding ideal. The group $G$ acts on
the
completion $\hat R=\limproj R/{\goth m}^n$.
Finally let $\widehat{R_0}$ be the ${\goth m}_0$-adic completion of
$R_0$.

The following easy lemma is certainly well known, but by of lack of reference
let me prove
the

\th LEMMA 2
\enonce With the previous notations, one has a canonical isomorphism
$$\hat R^G\iso {\hat R}_0$$.
\endth
\rem Remark
\endrem the lemma says that that passing to the invariants commutes with
completion.

Proof: one has to compare  $(\limproj R/{\goth m}^n)^G$ and $\limproj
R_0/{\goth m}_0^n$. Let $S$ be the coordinate ring of $G$ and denote by
$\sigma$ the morphism $$\sigma:\ R\ra S\otimes_k R$$
defining the action.  From [M], assertion (3) of theorem 1.2, the
ring of invariants of $R/{\goth m}^n$ is $ R_0/{\goth m}_0^n$. One therefore
has to prove that the
injection $$\limproj(R/{\goth m}^n)^G\mono(\limproj R/{\goth m}^n)^G$$
is onto. Let  $(r_n)\in \limproj R/{\goth m}^n$ an $G$-invariant sequence, that
is
$$\sigma(r_n)-1\otimes r_n\in S\otimes{\goth m}^{u_n}$$ for a sequence $u_n$
which goes to
$\infty$. After renormalization, one can assume that $u_n\geq n$ which implies
$\sigma(r_n)-1\otimes r_n\in S\otimes {\goth m}^n$ and therefore $(r_n)$
defines an element of
$\limproj(R/{\cal M}^n)^G$.\cqfd

Let $G_E$ denote the subgroup of
$G_q$ given by
$$G_E=\Ker\{\Aut(gr)\fhd{\det}{}G_m\}.$$
One can either prove by a direct calculation or by using Luna's result that
$G_E$ is
{\it reductive}. Let
$$A_E=k[\Ext^1_0(gr,gr)])^{G_E}$$
be the ring of polynomial maps on
$\Ext_0^1(gr,gr)$ invariant under $G_E$ (the group $G_E$ acts by functoriality
on both arguments
of $\Ext^1_0(gr,gr)$). Let ${\hat A_E}$ be the completion at the origin. Using
twice the previous
lemma and using Luna's theorem, one obtains

\th THEOREM 3
\enonce There is an isomorphism of complete local $k$-algebras
$$\hat\O_{\m, [E]}\iso {\hat A}_E.$$
\endth

\th COROLLARY
\enonce The local ring $\hat\O_{\m, [E]}$ depends only on the numerical
invariants of $X$ and
$gr(E)$.
\endth

\centerline{Suppose once for all that $E$ is non stable.}

One has of course the inequalities
$$1\leq\dim G_E\leq r^2-1\leqno(1)$$
with equality on the left hand side (resp. right hand side)  of (1) if
$G_E=G_m$ ({\bf Case 1})
(resp. $G_E=\Sl_r$ ({\bf Case 2})). Let's examine these 2 cases.

\section II. Case 1: $G_E=G_m$
\endsection

\ind In this case, the graduate $gr$ of $E$ is a direct sum
$$gr=gr_1\oplus gr_2$$
where $gr_i$ is stable of slope ${d\over r}$, rank $r_i\not=0$ and
$$gr_1\not\simeq gr_2.$$
Each element $(\alpha_1,\alpha_2)\in G_E(k)$ acts by mutiplication by
$\alpha_j.\alpha_i^{-1}$ on each factor $\Ext^1(gr_i,gr_j)$ of
$$\Ext^1(gr,gr)=\oplus \Ext^1(gr_i,gr_j).$$
Let
$$d_{i,j}=\dim \Ext^1(gr_i,gr_j)=\left\{
\matrix{
r_ir_j(g-1)&{\rm if}\ i\not=j\cr
r_ir_j(g-1)+1&{\rm if}\ i=j\cr}\right.$$
and $X_{i,j}^k,\ k=1,\ldots,d_{i,j}$ a basis of $\Ext^1_0(gr_i,gr_j)^\vee$.

The ring
$$A_E\subset k[\Ext^1(gr,gr)]=k[X_{i,j}^k,\ 1\leq i,j\leq n,\ 1\leq k\leq
d_{i,j}]$$
is the ring generated by
$\bl(\Ext^1(gr_1,gr_1)\oplus\Ext^1(gr_2,gr_2)\br)_0^\vee$ and the
products  $$<X_{1,2}^k.X_{2,1}^l,\ 1\leq k,l\leq d_{i,j}>.$$
Let $\S$ be the  the cone over the Segre variety
$$\P(\Ext^1(gr_1,gr_2))\times\P(\Ext^1(gr_2,gr_1))\subset
\P(\Ext^1(gr_1,gr_2)\otimes_k \Ext^1(gr_2,gr_1)).$$

\th PROPOSITION 1
\enonce There is an isomorphism
$$\Sp(A_E)\iso \bl(\Ext^1(gr_1,gr_1)\oplus\Ext^1(gr_2,gr_2)\br)_0\times \S.$$
\endth

Note that $\Sp(A_E)$ is a (homogeneous) cone. Using the theorem I.3 and the
previous
proposition, one therefore obtains the

\th PROPOSITION 2
\enonce With the previous notations,  the completion of $\m$ at $[E]$ is the
completion at the
origin of the cone  $\Sp(A_E)$.
\endth

Using furthermore that the multiplicity at the origin of the affine cone of a
projective variety
is just its degree, one gets the

\th COROLLARY 3
\enonce The completion of the tangent cone of $\m$ at $[E]$
is the completion at the origin of
$$\bl(\Ext^1(gr_1,gr_1)\oplus\Ext^1(gr_2,gr_2)\br)_0\times\S.$$
The Zariski tangent space is
$$\bl(\Ext^1(gr_1,gr_1)\oplus\Ext^1(gr_2,gr_2)\br)_0
\oplus\bl(\Ext^1(gr_1,gr_2)\otimes_k \Ext^1(gr_2,gr_1)\br).$$
Moreover the multiplicity of $\m$ at $[E]$ is
$$\mult_{[E]}(\m)=\pmatrix{2.d_{1,2}-2\cr d_{1,2}-1\cr}$$
(recall that $d_{1,2}=r_1.r_2(g-1)$).
\endth
\rem Remark 4
\endrem
Note that $\pmatrix{2.d_{1,2}-2\cr d_{1,2}-1\cr}=1$ if and only if $g=2$ and
$r_1=r_2=1$. Using that the singular locus of ${\rm Sing\ \m}$is closed in
$\m$, one obtains
easely in this way another proof of the fact that ${\rm Sing\ \m}$ is the non
stable locus, except
if $g=r=2$ and $d=0$.

\section III {\bf Case 2}: $G_E=\Sl_r$
\endsection

\ind In this case, $gr=L^{\oplus r}$ where $L^{\otimes r}=\O$. Using a
translation by $L^{-1}$
which induces an automorphism of  ${SU(r,\O)}$, one may assume $L=\O$. The ring
$A_E$ is the ring
of polynomial maps on $M_0(r)^g$ invariant under $\Sl_r$. This group  $\Sl_r$
acts diagonally by
conjugation on each factor $M_0(r)$, which is the space of traceless matrices
of size $r$. For
general $r$, Procesi [P] and Rasmyslev [Ras] have obtained the following
description of the first
2 syzigies of $A_E$:
\medskip
-Generators: for every sequence ${\underline i}=(i_1,\ldots,i_{N({\underline
i})})$ of integers of
$[1,\ldots,g]$, let $t_{\underline i}$ be the invariant polynomial map
$$t_{\underline i}:\ \left\{\matrix{
M_0(r)^g&\ra&k\cr
(X_1,\ldots,X_g)&\longmapsto&Tr(X_{i_1}.\ldots.X_{i_N})\cr
}.\right.$$ Then the $t_{\underline i}$ with $N({\underline i})\leq 2^g-1$ form
a system of
generators.
\medskip
-Relations between all the $t_{\underline i}'s$: Let $P_X$ be the
characteristic polynomial of the
general matrix $X$. The homogeneous polynomial $X\longmapsto Tr(X.P_X(X))$
gives by polarization
(namely by taking the
 total differential of order $g+1$) a multilinear map $F(H_1,\ldots,H_{g+1})$.
Then the relations
are generated by $tr_1,\ldots,tr_n$ and the relations $F(H_1,\ldots,H_{g+1})$
where the $H_i$'s
runs in the set of all possible monomials in the $X_i$.

Although this description is quite explicit, it looks difficult to obtain a
complete {\bf
finite} set of relations between the (finite) set of generators constructed
above.

As far as I know, the only case where such a finite descrition is available is
for $r=2$.
\footnote{(*)}{\sevenrm According to some experts of invariant theory, it is
more or less hopeless
to obtain such a finite description of $A_E$ in the general case}

\section III {\bf Case 2}: $G_E=\Sl_r$ and $r=2$
\endsection

In this case, $A_E$ can be described by using classical results of the
geometric invariant theory
of $SO_3(k)$. Following [LeB], §I.4, let me briefly explain this description.

For $X\in M_0(2)$, let $u(X)=\pmatrix{u_1(X),u_2(X),u_3(X)}\in k^3$ be defined
by the equality
$$X=\pmatrix{
u_1(X)&u_2(X)-\sqrt{-1}u_3(X)\cr
u_2+\sqrt{-1}u_3(X)&-u_1(X)\cr}.$$
By theorem 4.1 of [LeB] the isomorphism
$$\left\{\matrix{
M_0(2)&\ra&k^3\cr
X&\longmapsto&u(X)\cr}\right.$$
induces an identification of $A_E$ with the polynomial maps of $(k^3)^{\oplus
g}$
invariant under the canonical diagonal action of $SO_3(k)$.

Let $T_{i,j}$ be the invariant function corresponding to
$(u_1,\ldots,u_n)\longmapsto
(u_i.u_j)$ (scalar product), namely $T_{i,j}(X_1,\ldots,X_g)={1\over
2}Tr(X_iX_j)$.

Let $T_{i,j,k}$ be the invariant function corresponding to
$(u_1,\ldots,u_n)\longmapsto
u_i\wedge u_j\wedge u_k$ (the wedge product lives in
$\wedge{}^3k^3\build{=}_{}^{\rm can}k$),
namely $T_{i,j,k}(X_1,\ldots,X_g)=Tr(X_iX_jX_k)$.

With some abuse of notation, one can now use the results of H. Weyl [W],
theorem (2.9 A) and
(2.17 B) and it's sequel on page 77 which says the following:

-Generators for the invariants maps under ${O_3(k)}$: the set
$$<T_{i,j}>$$

-Relations between the generators:
the 4-minors of the $g\times g$-symmetric matrix
$$\pmatrix{&\vdots\cr
\ldots&T_{i,j}&\ldots\cr
&\vdots\cr
}$$
and the relations $T_{i,j}=T_{j,i}$. One recognizes the coordinate ring of
the (affine) cone $C$ of symmetric matrices of rank $\leq 3$. This scheme is
well understood: it
is integral and normal [easy], Cohen-Macaulay [H-R], its
multiplicity at the origin (or the degree of the projectivization $\P C$) is
known [H-T]...

-{\bf Generators} for the  ${SO_3(k)}$-invariants maps are: the $T_{i,j}$'s and
the
$T_{i,j,k}$'s.

-{\bf Relations}: the previous 4-minors and:

$$T_{i_1,i_2,i_3}.T_{j_1,j_2,j_3}=\det([T_{i_n,j_m}]_{1\leq n,m\leq 3})\leqno
(1)$$
$$T_{i_0,i_4}T_{i_1,i_2,i_3}-T_{i_1,i_4}T_{i_0,i_2,i_3}+T_{i_2,i_4}T_{i_0,i_1,i_3}-
T_{i_3,i_4}T_{i_0,i_1,i_2}\leqno (2)$$
and the relations given by the symmetry of $T_{i,j}$ and the skew symmetry of
$T_{i,j,k}$ in the
indices.

\rem Remark 1
\endrem
The relation (2) comes from the vanishing of the $4\times 4$-determinant
$\det([(u_{i_n}.v_{m})]_{0\leq n,m\leq 3})$ where $v_{0}=_{i_4}$ and
$v_{m}=e_m,\ m=1,2,3$ (the
vectors of the canonical basis of $k^3$).

\ind Let $\bar C$ be the tangent cone of $\Sp(A_E)$ at the origin. It is the
subscheme of
$$C\times \Sp k[T_{i,j,k}]$$
whose ideal is generated by
$$T_{i_1,i_2,i_3}.T_{j_1,j_2,j_3}\leqno (3)$$
and
$$T_{i_0,i_4}T_{i_1,i_2,i_3}-T_{i_1,i_4}T_{i_0,i_2,i_3}+T_{i_2,i_4}T_{i_0,i_1,i_3}-
T_{i_3,i_4}T_{i_0,i_1,i_2}\leqno (4)$$
and relations given by the skew symmetry of $T_{i,j,k}$ in the indices. The
ideal described
above is the ideal of initial forms of the ideals given by (1) and (2).

Let $k(C)$ be the function field of $C$ and $K$ its algebraic closure. Note
that, according to
(3), the ideal $I_{C/\bar C}$ of $C$ in $\bar C$ is nilpotent. This implies by
[F] example 4.3.4,
the formula $$\mult_0{\bar C}={\rm length}\O_{C,\bar C}.\mult_0(C).\leqno (5)$$
The next formula is clear
$${\rm length}\O_{C,\bar C}=1+\dim_{k(C)}I_{C/\bar C}\otimes_{\O_{\bar C}}
k(C)=1+
\dim_K I_{C/\bar C}\otimes_{\O_{\bar C}} K.\leqno (6)$$
\ind One therefore has to compute the
dimension over $K$ of the sub-vector space
$V_T$ of the dual space of $W=\oplus K.T_{i,j,k}$ of equations
given by (4) and the skew symmetry condition for the $T_{i,j,k}$. This vector
space is isomorphic
to  $I_{C/\bar C}\otimes_{\O_{\bar C}} K$.

\th LEMMA 2
\enonce The dimension of $V_T$ depends only on the conjugation class of $T$.
\endth

Proof: The symmetric matrix
$$T=[T_{i,j}]\in M_g(K)$$
acts on the dual vector space $V$ of $K^g$. Let $(e_i)_{1\leq
i\leq g}$ be  the canonical basis and  $(e_i^\vee)_{1\leq i\leq g}$ its dual
basis. The map

$$\left\{\matrix{
W&\ra&\wedge{}^3V\cr
T_{i,j,k}&\longmapsto&e_i\wedge e_j\wedge e_k\cr}\right.$$

identifies $W$ and $\wedge{}^3V^\vee$.
With this identification, the relations (4) become
$$-T(e_{i_4}^\vee)\int e_{i_0}\wedge e_{i_1}\wedge e_{i_2}\wedge e_{i_3}$$
and $\dim V_T$ is the corank of
$$\left\{\matrix{
V^\vee\otimes\wedge{}^4V&\ra&\wedge{}^3V\cr
x^\vee\otimes y&\longmapsto&T(x^\vee)\int y\cr}
\right.$$
This map depends only on the conjugation class of $T$.\cqfd

One can therefore assume that $T$ is diagonal of eigenvalues $\lambda_i$ with
$\lambda_i=0$ for $3<i\leq g$ and $\lambda_i\not=0$ if $i\leq 3$.

Let us prove this simple lemma

\th LEMMA 3
\enonce The dimension of $V_T$ is ${g(g-1)(g-2)\over 6}-{(g-3)(g-4)(g-5)\over
6}$ if $g\geq 3$ and
$0$ if $g\leq 2$.
\endth

Proof: If $g\leq 2$,  the vector space $\wedge{}^3V$ is zero and so is $V_T$.
Suppose $g\geq 3$.
Consider an equation
$$T_{i_0,i_4}T_{i_1,i_2,i_3}-T_{i_1,i_4}T_{i_0,i_2,i_3}+T_{i_2,i_4}T_{i_0,i_1,i_3}-
T_{i_3,i_4}T_{i_0,i_1,i_2}\leqno (4)$$ defining $V_T$. If $i_4\not\in\{1,2,3\}$
or
$i_4\not\in\{i_0,i_1,i_2,i_3\}$ then the equation (4) is trivial.

Let me suppose that $i_4\in\{1,2,3\}$ and for instance that $i_4=i_0$. If
$i_4\in\{i_1,i_2,i_3\}$,
the equation is just a consequence of the skew symmetry of $T_{i_1,i_2,i_3}$.
To get a new
relation, one has therefore further to suppose further that
$i_4\not\in\{i_1,i_2,i_3\}$ and the
equation become $$\lambda_{i_0}.T_{i_1,i_2,i_3}=0.$$ Of course, the other cases
are obtained by
symmetry. One has proved the following:

the equations (4) are non trivial if and only if
$$\{i_0,i_1,i_2,i_3\}=\{i,j,k\}\sqcup\{i_4\}\ {\rm and\ }i_4\not\in\{i,j,k\}.$$
In this case the relation (4) becomes
$$T_{i,j,k}=0,$$
or equivalently
$$T_{i,j,k}=0\ {\rm if\ }\{i,j,k\}\cap\{1,2,3\}\not=\emptyset.$$ In particular
this corank is
$$\dim\wedge{}^3V-A^3_{g-3}={g(g-1)(g-2)\over 6}-{(g-3)(g-4)(g-5)\over
6}.$$\cqfd

The degree $d_g^r$ of the locus
\footnote{(*)}{\sevenrm In [H-T], this locus is endowed with the reduced scheme
structure. But it
is known in full generality that the natural scheme structure given by the
vanishing of the
$(g-r+1)$-minors is Cohen-Macaulay [J] and generically reduced [easy] and
therefore reduced. In
our case ($r=g-3$), this reduceness is obvious, because $C$ is a ring of
invaraiants.}
of $g\times g$-symmetric
matrices of corank $\geq r$ is computed in [H-T], proposition 12.b:
$$d_g^r=\deg\P C=\prod_{\alpha =0}^{r-1}
{\pmatrix{g+\alpha \cr r-\alpha \cr}\over\pmatrix{2\alpha +1\cr\alpha \cr}}.$$

Using the formulas (5) and (6) and the lemma 3, one obtains the

\th THEOREM 4
\enonce The multiplicity of $[\O\oplus\O]$ in $SU(2,\O)$ is
$$\Bl(1+{g(g-1)(g-2)\over 6}-{(g-3)(g-4)(g-5)\over 6}\Br).d_g^{g-3}$$
if $g\geq 3$ and $1$ if $g=2$.
\endth

\rem Remarks 5
\endrem

1.-The preceding discussion gives a precise description of the tangent cone
$\bar C$ of
$SU(2,\O)$. In particular, the Zariski tangent space at the trivial bundle is
$$T_{[\O\oplus\O]}SU(2,\O)={\rm Sym}^2V\oplus\wedge{}^3V.$$

2.-One recovers the smoothness of $SU(2,\O)$ if $g=2$.
\medskip
\ind In the case of a a non hyperelliptic genus $3$  curve, the generalized
$\Theta$ divisor
embeds  $SU(2,\O)$ as the  Cobble quartic $SU(2,\O)$ in $\P
H^0(J^2,2.\Theta_{J^2})$ (see [N-R2]
and [D-O] pages 184-185). Let $q_X$ be an equation of the Coble quartic.

\def\A{{\bf A}}
\th THEOREM 6
\enonce Suppose that $X$ is a non hyperelliptic genus $3$  curve.

(i) The local equation of
the Coble quartic at the trivial bundle is
$$T^2=\det([T_{i,j}]_{1\leq i,j\leq 3})$$ in the
affine space ${\bf A}^7$ with coordinates $T,T_{i,j}$ with $T_{i,j}=T_{j,i}$.

(ii)The local equation of the Coble quartic at $gr=gr_1\oplus gr_2$ with
$gr_1\not=gr_2$ is a
rank 4 quadric in $\A^7$.

(iii) The ideal of the Kummer $K(X)$ variety of $J(X)$ in $\P
H^0(J^2,2.\Theta_{J^2})$ is generated by the 8 cubic equations which are the
partials derivatives
of the Coble quartic. \endth

Proof: The first 2 points are clear from proposition II.2 and (1), (2) . Let me
prove
({\it iii}). Let  $\K$ the scheme  defined by the partials of $q_X$. The Kummer
variety $K(X)$ is
the reduced scheme of $\K$. It is therefore enough to prove that the completion
of $\K$
at each non stable class $[gr]$ of $K(X)\subset\K$ is {\it reduced}.

Because of the invariance of the Coble quartic under the Theta group of
$2.\Theta_{J^2}$, there
are 2 cases: either $gr$ is trivial, or $gr=gr_1\oplus gr_2$ with
$gr_1\not=gr_2$. In the
first case, by ({\it i}), the equations in $k[[T,T_{i,j}]]$ of the completion
of $\K$ at
$[\O\oplus\O]$ are $T$ and the $2\times 2$-minors of $[T_{i,j}]_{1\leq i,j\leq
3}$. It is
precisely the (completion at the origin) of the cone over the Veronese surface
in $\P^5$ (with
homogenous coordinates $T_{i,j}$) and $K$ is therefore reduced. The second case
is even
simpler, $\K$ being (the completion of) a $3$-plane in $\A^7$ (the tangent
space of
$K(X)$).\cqfd

\def\A{{{\cal A}_{gr}}}\def\n{{\underline n}}\def\m{{\underline m}}
\def\p{{\goth p}}\def\mod{\mathop{\rm\ \ mod\ \ }\nolimits}\def\M{{\cal M}}
\section III The case $SU(3,\O)$ for of a genus $2$ curve
\endsection

\ind Suppose in this section that $X$ has genus $2$ and let $\M=SU(3,\O)$ be
the moduli space of
rank 3 semi-stable vector bundles on $X$ with trivial determinant. Consider a
non stable point of
$\M$ defining by the graded object $gr$ of a semi-stable bundle.

\ind The case $G_{gr}=G_m$ has been treated in section
II: in this case, the completion of $\M$ at $gr$ is the completion at the
origin of a rank 4
quadric in ${\bf A}^9$.
\medskip
\ind Suppose now that $G_{gr}=G_m\times G_m$ which means
$$gr=gr_1\oplus gr_2\oplus gr_3\ {\rm with\ }gr_i\not=gr_j{\ \rm for\
}i\not=j{\rm\ and\ }\deg
(gr_i)=0.$$

Let $X_{i,j}$ be a basis of $\Ext^1(gr_i,gr_j)^\vee$ for $i\not=j$. Let
$$\A=k[X_{i,j}]^{G_{gr}}$$
be the ring of invariants of $k[X_{i,j}]$ under $G_{gr}$ with the action
defined by the
following rule
$${\underline\alpha}.X^\n=\prod\bl({\alpha _i\over\alpha
_j}\br)^{n_{i,j}}X^\n$$
for
$(\alpha_1,\alpha _2,\alpha 3)\in G_{gr}(k)=(k^*)^3/k^*$. The following
equality is easy
$$A_{gr}=k[\bl(\oplus\Ext^1(gr_i,gr_i)\br)_0]\otimes\A.$$

A polynomial $\sum p_{\underline n}X^\n$ is in $\A$
if and only if
$$\sum_j n_{j,i}=\sum_j n_{i,j}\ {\rm for}\ i=1,2,3\leqno (1)$$
if $p_{\underline n}\not=0$. Therefore, $\A$ is generated by the monomials
$$X^\n\ {\rm such\ that\ }\n{\rm\ satisfies\ }(1).$$

\th LEMMA 1
\enonce The ring $\A$ is generated by
$$X_{3,2}X_{2,1}X_{1,3},\ X_{1,2}X_{2,3}X_{3,1}\ {\rm and\ }X_{i,j}X_{j,i},\
i<j.$$
\endth

Proof: put $\delta_{i,j}=n_{i,j}-n_{j,i}$ and let $\n$ be a multi-index
satisfying (1). The
relations (1) become
$$\delta_{1,2}+\delta_{1,3}=0,\ \delta_{1,2}=\delta_{2,3},\
\delta_{1,3}+\delta_{2,3}=0.$$

If $\delta^+=\delta_{1,2}\geq 0$, we write
$${\underline
n}=(n_{2,1}+\delta^+,n_{2,1},n_{1,3},n_{1,3}+\delta^+,n_{3,2}+\delta^+,n_{3,2})$$
and use the monomial
$$X_{1,2}X_{2,3}X_{3,1}$$
corresponding to
$${\underline n}_0=(1,0,0,1,1,0)$$
to
write ${\underline n}={\underline m}+\delta^+.{\underline n}_0$. This allows us
to write

$$X^\n=(X_{1,2}X_{2,3}X_{3,1})^{\delta^+}\prod_{i<j}(X_{i,j}.X_{j,i})^{m_{i,j}}$$
with $m_{i,j}\geq 0$. In the same way, when $\delta^-=-\delta_{1,2}>0$, one has
an equality
$$X^\n=(X_{3,2}X_{2,1}X_{1,3})^{\delta^-}\prod_{i<j}(X_{i,j}.X_{j,i})^{m_{i,j}}$$
with $m_{i,j}\geq 0$.\cqfd

For $i\in I=\{1,2,3\}$, put
$$x_i=X_{j,k},\ x_{i+3}=X_{k,j},\ \zeta_i=x_ix_{i+3}$$
with $I=\{i,j,k\}$ and $j<k$.Let  $\zeta_4=X_{3,2}X_{2,1}X_{1,3}$ and
$\zeta_5=X_{1,2}X_{2,3}X_{3,1}$. There is an equality
$$\zeta_4\zeta_5=\zeta_1\zeta_2\zeta_3.$$

\th PROPOSITION 2
\enonce The natural morphism
$$f:\ \left\{
\matrix{
k[X_i]&\ra&\A\cr
X_i&\longmapsto&\zeta_i\cr
}\right.$$
gives an
isomorphism
$$k[X_i]/(X_4X_5-X_1X_2X_3)\iso\A.$$
\endth

Proof: Let $\p$ be the (prime) ideal generated by $X_4X_5-X_1X_2X_3$.
Let
$$P=\sum\alpha_\n X^\n$$
be an element of $\Ker(f)$. Then, one finds by simple expansion
$$f(P)=\sum\alpha_\n \prod
x^{\phi(\n)}=\sum_{\m}x^\m\sum_{\phi(\n)=\m}\alpha_\n$$
with
$$\phi(\n)=(n_1+n_5,n_2+n_4,n_3+n_5,n_1+n_4,n_2+n_5,n_3+n_5)$$
which implies
$$\sum_{\phi(\n)=\m}\alpha_\n=0.\leqno (2)$$
(Here $X^\n=\prod_i X_i^{n_i}$ and $x^\m=\prod_i x_i^{m_i}$
$\n=(n_i)_{1\leq i\leq 5}$ and $\m=(m_i)_{1\leq i\leq 6}$ are multi-indices).

The kernel of $\phi$ is generated by $(1,1,1,-1,-1)$:

if $\phi(\m')=\phi(\m)$, there exists $\alpha\in\Z_+$ such that
$$\pm\alpha (1,1,1,0,0)+\m'=\pm\alpha (0,0,0,1,1)+\m.$$
In particular, one has the congruence
$$(X_4.X_5)^\alpha.X^{\m'}\equiv (X_4.X_5)^\alpha.X^\m\mod\p.\leqno (3)$$
According to (2) and (3),
we get the existence of a positive integer $a$ such that $$(X_4.X_5)^a.P\equiv
0\mod\p.$$
Since the ideal $\p$ is prime and $(X_4.X_5)\not\in\p$ and therefore
$P\in\p$.\cqfd

\th COROLLARY 3
\enonce The completion of $\M$ at a point $[gr]$ satisfying $G_{gr}=G_m\times
G_m$
is the completion at the origin of
$$\bl(\oplus\Ext^1(gr_i,gr_i)\br)_0\times\Sp(k[X_i]/(X_4X_5-X_1X_2X_3)).$$
Its tangent cone is a rank 2 quadric in the Zariski tangent space
$T_{[E]}\M={\bf
A}^9$.
\endth

\ind In particular, there exists a family $\bf E$ of semi-stable bundles of
trivial
determinant over a germ of curve such that:

(i) The group $G_{{\bf E}_\eta}$ of the generic bundle ${\bf E}_\eta$ is
$G_m\otimes_k k(\eta)$.

(ii) The group of $G_{{\bf E}_s}$ the special bundle ${\bf E}_s$ is $G_m\times
G_m$.

(iii) The multiplicity of $\M$ at $[{\bf E}_{\eta}]$ and $[{\bf E}_s]$ are the
same.

This shows that 2 points of $\M$ can have the same multiplicity without having
the same group
of automorphisms.

\ind When the $gr$ has 3 summands for which at least 2  are isomorphic, or
equivalently if
$G_{gr}$ is not a torus, the calculations are very intricate (but seem to be
possible). In fact,
one can in spite of this obtain the following

\th PROPOSITION 4
\enonce The tangent cone at each non stable point $gr$
such that $G_{gr}$ is not a torus is a quadric in ${\bf A}^9$ of rank $\leq 2$.
\endth

Proof: thanks to the corollary 1.7.4 of [Ray], for {\it every} semi-stable
vector
bundle $E$ of rank 3 and determinant $\O$, the determinantal locus $\Theta_E$
in $\Pic^{1}(X)$
$$\Theta_E=\{L\in\Pic^{1}(X)\ {\rm such\ that}\ H^0(X,E\otimes  L)\not=0\}$$
is a {\it divisor} in
$\mid\!\! \O(3\Theta_{J^1})\!\!\mid$ where $\Theta_{J^1}$ is the canonical
theta divisor on
$\Pic^{1}(X)$. The Picard group of $\M$ is cyclic with ample generator
$\O(\Theta)$ [D-N]. By
[B-N-R], the inverse image of $\O(1)$ by the morphism
$$\pi:\ \left\{\matrix{
\M&\ra&\mid\!\! \O(3\Theta_{J^1})\!\!\mid\cr
[E]&\longmapsto&\Theta_E\cr
}\right.$$
is $\O(\Theta)$  and $\pi^*$ is
a (canonical) isomorphism $$\mid\!\! \Theta\!\!\mid^\vee\iso\mid\!\!
3.\Theta_{J^1}\!\!\mid.$$
(In particular, $\mid\!\! \Theta\!\!\mid$ has no base point in this case!)

Using this isomorphism, $\pi$ becomes the morphism given
by the complete linear system $\mid\!\! \Theta\!\!\mid$. There are various ways
to prove this
simple lemma, but the following one can be generalized for the higher rank
case.

\th LEMMA 5\ \footnote{(*)}{\sevenrm This fact, which is due to {\pc I.} {\pc
DOLGACEV},
was pointed out to me by {\pc A.} {\pc BEAUVILLE}.}
\enonce The morphism $\pi$ is finite of degree $2$ over $\P^8$.
\endth

Proof of the lemma: using the isomorphism $\pi^*\O(1)\iso\O(\Theta)$, we get
that $\pi$ is
 finite of degree $c_1(\Theta)^8$ onto $\P^8$. One therefore has to compute the
degree of $\pi$. Although there exists a general beautiful formula due to
Witten to evaluate the
volume $${\displaystyle c_1(\Theta)^{\dim SU(r,\O)}\over\displaystyle (\dim
SU(r,\O))!},$$ we give
a simplest (and elementary) method to get this volume for $\M$. One has to
prove that the leading
term of the Hilbert polynomial
$$n\longmapsto P(n)=\chi(X,\Theta^n)$$
is ${2\over 8!}$. The canonical divisor of $\M$ si $\Theta^{-6}$ [D-N]. Serre
duality implies
therefore the symmetry $$P(n)=P(-6-n).\leqno (4)$$
The Grauert-Reimenschneider vanishing theorem (recall that $\M$ has rational
singularities) gives
the equality $$P(n)=\dim H^0(\M,\Theta^n)\ {\rm for\ }n\geq -5.$$
One therefore obtains the values
$$P(n)=0\ {\rm for\ }n=-5,\ldots,-1,\ P(0)=1{\rm\ and\ }P(1)=9.\leqno(5)$$
By (4) and (5), one obtains
$$P(X)=\lambda(X+5)(X+2)(X+3)^2(X+2)(X+1)(X-\alpha)(X+6+\alpha).$$
The equalities $P(0)=1{\rm\ and\ }P(1)=9$ imply
$$\alpha =-3\pm\sqrt{-47}{\rm\ and\ }\lambda={2\over 8!}.$$\cqfd

One has proved that the morphism $\pi$ is finite of degree 2 onto
$\mid\!\!\Theta\!\!\mid^\vee=\P^8$. Since $\M$ is Cohen-Macaulay and $\P^8$
smooth, this double
covering is flat (E.G.A. IV.15.4.2) and is given locally by an equation
$$t^2=f({\underline x})$$
where $\underline x$ are local coordinates on $\P^8$. This implies that the
multiplicity of each
point of $\M$ is $\leq 2$.

Take a point $[gr]\in\M$ such that $G_{gr}$ is not a torus: it is a non smooth
point of $\M$,
therefore the tangent cone is a quadric (the initial term of $f$ is not
linear). But, such a point
is a specialization of a point $gr_\eta$ such that $G_{gr_\eta}=G_m\times G_m$:
by the (obvious)
semi-continuity of the rank of the quadric cone of $[gr]$, the inequality ${\rm
rank}\leq 2$ follows
from the corollary 3.\cqfd

\rem Remark 6
\endrem
There exists by the way a nice $3$-dimensional family of sextics in $\P^8=\P
H^0\bl(J^1,\O(3.\Theta_{J^1})\br)$ (given by the ramification of $\pi$), which
could be called the
family of Dolgacev-Coble sextics associated to a genus 2 curve. It would be
interesting to know if
these sextics share the same kind of properties as the Coble quartics.

%\end

%\input utilitaire_locale
\section{IV. Multiplicity of the theta divisor (rank 2 case)}
\endsection

\ind Recall that there exists a (Cartier) divisor $\Theta$ on $SU(2,\omega_X)$
which is
characterized by the following universal property [D-N]:
let $S$ be a $k$-scheme and ${\bf E}$ a family of semi-stable vector bundles
over
$X_S=X\times_k S$ of determinant $\omega_X$. Let $\pi:\ S\ra SU_X(2,\omega_X)$
be the classifying
map corresponding to $\bf E$. Then, one has the equality
$$\pi^*(\Theta)=\div (\det Rp_*{\bf E})^\vee.$$

By construction, the geometric points of $\Theta$ are the classes $[E]\in
SU(2,\omega_X)$ such that
$H^0(X,E)\not=0$.

\rem Remark
\endrem The vanishing of $H^0(X,E)$ depends only on the S-equivalence class
$[E]$ of $E$.

Let $E$ be a semi-stable vector bundle of determinant $\omega_X$. Recall ([La],
theorem III.3)
that for $E$ stable
$$\mult_{[E]}\Theta=\dim H^0(X,E)$$
and that the tangent cone is defined in $\Ext_0^1(E,E)$ by the ideal of the
determinant of linear
forms defined by the cup-product
$$H^0(X,E)\otimes \Ext^1_0(E,E)\ra H^1(X,E).$$

\ind These facts can be generalized formally (using the universal property of
$\Theta$) as
follows. Let  $[E]\in SU(2,\omega_X)$ a non stable point of $\Theta$ of graded
object
$gr=gr_1\oplus gr_2$ and let $$h={1\over 2}\dim H^0(X,gr)=\dim H^0(X,gr_1)=\dim
H^0(X,gr_2)$$
(note  that $gr_1\otimes gr_2=\omega_X$ which
implies by Serre duality  and Riemann-Roch the equality
$$\dim H^0(X,gr_2)=\dim H^1(X,gr_2)=\dim H^0(X,gr_1)).\leqno (1)$$

With the notations of §I, let $V=\Sp k[[\Ext^1_0(gr,gr)]]$ be a (formal)
\'etale slice of $\cal Q$
at $gr$ and  $$\pi:\ V\ra V/\Aut(gr)\fhd{\hbox{\'etale}}{} SU(2,\omega_X)$$ the
canonical
morphism.  Then, the induced map
$$\pi^*(\Theta)/\Aut(gr)\ra\Theta$$
is \'etale. The tangent cone of $\pi^*(\Theta)$ is given by the determinant
$$d_{gr}\in {\rm Sym}^{2h}\Ext^1_0(gr,gr)^\vee$$
defined (up to a non zero scalar) by the cup-product
$$H^0(X,gr)\otimes \Ext^1_0(gr,gr)\ra H^1(X,gr)$$
In particular, a point $e\in\Ext^1_0(gr,gr)$ is in the tangent cone of
$\pi^*(\Theta)$ at $[E]$ if
and only if the cup-product
$$\cup e:\ H^0(X,gr)\ra H^1(X,gr)$$ is non onto.

\th PROPOSITION 1
\enonce Assume that $gr_1\not=gr_2$. Then, the multiplicity of $\Theta$ at
$[gr]$ is
$$\mult_{[gr]}\Theta={1\over 2}\dim H^0(X,gr).\mult_{[gr]}SU(2,\O).$$
\endth

Proof: with the notation of the second section, the completion of
$SU(2,\omega_X)$ at $[gr]$ is
the completion at the origin of
$$\bl(\Ext^1(gr_1,gr_1)\oplus\Ext^1(gr_2,gr_2)\br)_0\times \S$$
where $\S$ is the  the cone over the Segre variety
$$\P(\Ext^1(gr_1,gr_2))\times\P(\Ext^1(gr_2,gr_1))\subset
\P(\Ext^1(gr_1,gr_2)\otimes_k \Ext^1(gr_2,gr_1)).$$
Fix coordinates
$${\underline X}=(X_{i,j}^k){\rm\ on\ }\Ext^1(gr_i,gr_j){\rm\ for\ }i\not=j\
{\rm and\ }
(\underline Y)=Y_{i}^k{\rm\ on\ }
\Ext^1(gr_i,gr_i).$$

The equation $F$ of $\pi^*(\Theta)$ is of the form
$$F=d_{gr}+G_{2h+1}$$
where $G_{2h+1}$ vanishes at the origin with order $\geq 2h+1$. The polynomials
$d_{gr}$ and
$G_{2h+1}$ are $G_E$ invariant and therefore (see section 2) can be writen in
terms of
$\underline  Y$ and $z^{k,l}=X_{1,2}^k.X_{2,1}^l$. Let me decompose $d_{gr}$ as
$$d_{gr}=\sum_{i=0}^{2h} Q_{i}({\underline X})P_{2r-i}({\underline Y})$$
where the degree of $Q_i$ (resp. $P_{2r-i}$) is $i$ (resp. $2r-i$) and $P_0=1$.
Using the
degrees, one finds the following properties:

-If $i$ is odd, then $Q_{i}({\underline X})P_{2r-i}({\underline Y})=0$.

-If $P_{2r-2i}\not=0$, then $Q_{2i}$ is invariant and therefore
$$Q_{2i}({\underline X})=R_i({\underline z})$$
with $R_i$ is a polynomial in $\underline z$ of degree $i$ which is defined up
to the ideal of the
Segre cone $\S$.

It follows that the equation of $\pi^*(\Theta)$ can therefore be written as
$$R_h({\underline z})+S({\underline z},{\underline Y})\leqno (2)$$
where $S$ vanishes at the origin with order $\geq h+1$ at the origin.

\th LEMMA 2
\enonce The polynomial $Q_{2h}(\underline X)$ is non zero.
\endth

Proof of the lemma: according to the previous discussion, one just has to prove
the existence of
$$e\in \Ext^1(gr_1,gr_2)\oplus\Ext^1(gr_2,gr_1)\subset\Ext^1_0(E,E)$$
such that the cup product
$$\cup e:\ H^0(X,gr)\ra H^1(X,gr)$$
is onto. By symmetry, one only has to prove the existence of $e_1\in
\Ext^1(gr_1,gr_2)$
such that the cup product
$$\cup e_1:\ H^0(X,gr_1)\ra H^1(X,gr_2)$$
is onto. This is classical (see [La], lemma II.8): let $\Gamma$ be the variety
$$\Gamma=\{(k.s,k.e)\in\P H^0(X,gr_1)\times\P\Ext^1(gr_1,gr_2)\ {\rm such\
that}\ s\cup e=0\}$$
and $p$ (resp. $q$) the first (resp. second) projection. Let $0\not=s\in
H^0(X,gr_1)$ and
$D=\div(s)$ its
zero divisor. The canonical surjection
$$\cup s:\ {\cal H}om(gr_1,gr_2)\epi gr_2(-D)$$
gives a surjection
$$\cup s:\ \Ext^1(gr_1,gr_2)\epi H^1(X,gr_2).\leqno (3)$$
By (3) the dimension of $p^{-1}(k.s)$ is
$$\dim p^{-1}(k.s)=\dim\P\Ext^1(gr_1,gr_2)-\P\dim H^0(X,gr_1)-1$$

Therefore $$\dim\Gamma=\dim\P\Ext^1(gr_1,gr_2)-1$$ and $q(\Gamma)\not=\Gamma.$
\cqfd
The polynomial $R_h$ can be thought of as an element of
$$H^0(\P\Ext^1(gr_1,gr_2)\times\P\Ext^1(gr_2,gr_1),\O(h,h))$$ which by the
lemma 2 is non zero.
According to (2), the completion of the tangent cone is therefore the
hypersurface of
$$\bl(\Ext^1(gr_1,gr_1)\oplus\Ext^1(gr_2,gr_2)\br)_0\times \S$$
given by $R_h$. The proposition follows.\cqfd

\section{REFERENCES}
\endsection
\medskip
[B-N-R] A. {\pc BEAUVILLE}, M. S. {\pc NARASIMHAN}, M. S. {\pc RAMANAN} , {\it
Spectral curves
and the generalized theta divisor}, J. reine angew. Math., {\bf 398}, 169-179
(1989).
\medskip
[B] J. F. {\pc BOUTOT}, {\it Singularités rationelles et quotients par les
groupes réductifs},
Inv. Math;, {\bf 88}, 65-68 ((1987).
\medskip
[D-N] J. M. {\pc DREZET}, M. S. {\pc NARASIMHAN}, {\it Groupe de Picard des
vari\'et\'es de
modules de fibr\'es semi-stables sur les courbes alg\'ebriques}, Invent. math.,
{\bf 97}, 53-94
(1989).
\medskip
[D-O] I.{\pc DOLGACEV}, D. {\pc ORTLAND}, {\it Point sets in projective spaces
and theta functions}, Ast\'erisque, {\bf 165}, (1988).
\medskip
[F] W. {\pc FULTON} {\it Intersection theory}, Ergeb. der Math. und ihrer
Grenzgebiete,  3. Folge
Band 2,  Springer-Verlag, Berlin Heidelberg New-York Tokyo (1984).
\medskip
[H-R] M. {\pc HOCHSTER}, J. {\pc ROBERTS}, {\it Rings of invariants of
reductive
groups acting on regular local rings are Cohen-Macaulay}, Adv. in  Math. , {\bf
13}, 115-175
(1974).
\medskip  [H-T] J. {\pc HARRIS}, L. W. {\pc TU}, {\it On symmetric and skew
symmetric
determinantal varieties}, Topology, {\bf 23}, 71-84 (1984).
\medskip
[Ha] R. {\pc HARTSHORNE}, {\it Algebraic Geometry}, GTM 52, Spinger Verlag,
Berlin-Heidelberg-New York, (1977).
\medskip
[J] T. {\pc JOZEFIAK}, {\it Ideals generated by minors of a symmetric matrix},
Comment. Math.
Helvetici, {\bf 53}, 595-607 (1978).
\medskip
[La] Y. {\pc LASZLO}, {\it Un th\'eor\`eme
de Riemann pour les diviseurs Th\^eta g\'en\'eralis\'es sur les espaces de
modules de fibr\'es
stables sur une courbe}, Duke Math. Journ., {\bf 64}, 333-347 (1991).
\medskip
[LeB] L. {\pc LE \pc BRUYN}, {\it Trace rings of generic 2 by 2 matrices}, Mem.
A.M.S., {\bf 66},
(1987).
\medskip
[Lu] D. {\pc LUNA}, {\it Slices \'etales}, M\'em. Soc. math. France, {\bf 33},
81-105 (1973).
\medskip
[N-R1] M. S. {\pc NARASIMHAN}, M. S. {\pc RAMANAN} , {\it Moduli of vector
bundles on a
compact Riemann surface}, Ann. of Math., {\bf 89}, 14-51 (1969).
\medskip
[N-R2] M. S. {\pc NARASIMHAN},M. S. {\pc RAMANAN} , {\it $2\theta$-linear
system on abelian
varieties, in} Vector bundles and algebraic varieties, Oxford University Press,
415-427 (1987).
\medskip
[P] C. {\pc PROCESI} {\it The invariant theory of $n$ by $n$ matrices}, Adv. in
Math. {\bf 19},
306-381.
\medskip
[Ras] J. {\pc RASMYSLEV}, {\it Trace identities of full matrix algebras over a
field of
characteristic zero}, Izv. Akad. Nauk. U.S.S.R., {\bf 8}, 727-760 (1974).
\medskip
[Ray] M. {\pc RAYNAUD}, {\it Sections des fibr\'es vectoriels sur une courbe},
Bull.
Soc. math. France, {\bf 110}, 103-125 (1982).
\medskip
[S] C. S. {\pc SESHADRI}, {\it Fibr\'es vectoriels sur les courbes
alg\'ebriques} (r\'edig\'e par
J. M. {\pc DREZET}), Ast\'erisque, {\bf 96}, (1982).
\medskip
[W] H. {\pc WEYL}, {\it The classical groups},
Princeton Univ. Press, Princeton, N.J., (1946).
\footnote{}{\sevenrm Yves LASZLO  Universit\'e Bordeaux 1 Math\'ematiques 351,
cours de la Lib\'eration 33405 Talence Cedex France and Universit\'e Paris-Sud
U.R.A. 752
Math\'ematiques B\^atiment 425 91405 Orsay Cedex France}

\end